\documentclass[12pt]{article}
\usepackage{amsmath}
\usepackage{amssymb}

\def\qed{\hfill$\Box$\par\medskip\par\relax}

\numberwithin{equation}{section}
\newcommand{\eps}{\varepsilon}
\newcommand{\Z}{{\mathbb Z}}
\newcommand{\N}{{\mathbb N}}
\newcommand{\V}{{\mathcal V}}
\newcommand{\K}{{\mathcal K}}
\newcommand{\B}{{\mathcal B}}
\newcommand{\M}{{\mathcal M}}

\newcommand{\R}{{\mathbb R}}

\newcommand{\ZZ}{{\mathcal Z}}

\renewcommand{\phi}{\varphi}
\newcommand{\A}{{\mathfrak A}}

\newcommand{\sgn}{\mathop{\rm sgn}\nolimits}

\newcommand{\IE}{{\mathbb E}}
\newcommand{\IP}{{\mathbb P}}
\newcommand{\Po}{{\mathtt P}_{\!\omega}}
\newcommand{\Eo}{{\mathtt E}_{\omega}}
\newcommand{\PoA}[1]{{\mathtt P}_{\omega|#1}}

\allowdisplaybreaks

\newtheorem{theo}{Theorem}[section]

\newtheorem{df}[theo]{Definition}

\newtheorem{rmk}[theo]{Remark}

\title{Survival of branching random walks in random environment}

\author{Nina Gantert$^{1}$ \and Sebastian M\"{u}ller$^{2}$ \and
 Serguei~Popov$^{3}$ \and Marina Vachkovskaia$^{3}$}

\begin{document}

\maketitle

{\footnotesize
\noindent $^{~1}$CeNos Center for Nonlinear Science and Institut f\"{u}r Mathematische Statistik,
Fachbe\-reich Mathematik und Informatik,
Ein\-stein\-strasse 62, 48149 M\"{u}nster,
Germany

\noindent e-mail: \texttt{gantert@math.uni-muenster.de}

\noindent $^{~2}$  Institut f\"{u}r Mathematische Strukturtheorie,
Technische Universit\"{a}t Graz, Steyrergasse 30, 8010 Graz, Austria

\noindent e-mail: \texttt{mueller@tugraz.at}, \quad \noindent url:
\texttt{http://www.math.tugraz.at/$\sim$mueller}

% \noindent $^{~3}$Instituto de Matem{\'a}tica e Estat{\'\i}stica,
% Universidade de S{\~a}o Paulo, rua do Mat{\~a}o 1010, CEP
% 05508--090, S{\~a}o Paulo, SP, Brazil
% 
% \noindent e-mail: \texttt{popov@ime.usp.br}, \quad
% \noindent url: \texttt{http://www.ime.usp.br/$\sim$popov}

\noindent $^{~3}$ Department of Statistics,
Institute of Mathematics, Statistics and Scientific Computation,
University of Campinas--UNICAMP,
P.O.\ Box 6065, CEP 13083--970, Campinas, SP, Brazil

\noindent e-mail: \texttt{$\{$popov,marinav$\}$@ime.unicamp.br}, \newline
\noindent url: \texttt{http://www.ime.usp.br/$\sim$$\{$popov,marinav$\}$}

}

\begin{abstract}
We study survival of nearest-neighbour branching random walks
in random environment (BRWRE) on $\Z$. A priori there are three
different regimes of survival: global survival, local survival,
and strong local survival. We show that local and strong local
survival regimes coincide for BRWRE and that they can be characterized with
the spectral radius of the first moment matrix of the process.
These results are  generalizations of the classification of BRWRE
in recurrent and transient regimes.
Our main result is a
characterization of global survival that is given in terms of
Lyapunov exponents of an infinite product of i.i.d.\ $2\times 2$ random matrices. 
\end{abstract}

\section{Introduction}
\label{s_intro}
% Let us first describe the process that evolves on~$\Z$ in discrete
% time. 
The branching random walk in random environment (BRWRE)
starts with one particle in the origin. This particle splits up in
several other particles at positions $\{-1,0,1\}$ according to some offspring
distribution. Now the process is defined
inductively, at each moment each particle at~$x$, 
independently of the other particles and the history of the process, splits up in
particles at $\{x-1,x,x+1\}$ according to some offspring
distribution that may depend on~$x$. The collection of the offspring
distributions itself is chosen randomly before starting the
process and then is kept fixed during the evolution of the
process. 
The difference with the model in \cite{BCGH2, GH2, GH3} is that we start the process with one particle and not with infinitely many. In contrast to previous papers on the topic, cf.\
\cite{CMP, CP, CP07, MP1, MP2, Mueller}, we allow the process to die out. A priori there
are three different types of survival: global survival, local survival, and strong
%N
local survival. Global survival means that with positive probability the
total number of particles is always positive. Local survival means
that with positive probability every site is visited
infinitely often. If the two latter probabilities are equal and
%N
positive (so that, conditioned on always having a positive number of particles, every site
is visited infinitely often a.s.)
we say there is strong local survival. Our first
result (Theorem~\ref{s_loc-ext})
says that, in fact, local and strong local survival coincide and
do not depend on the realization of the environment. 
%N
This is a generalization of the
classification of BRWRE in recurrent and transient regimes, 
compare with~\cite{CP} and~\cite{Mueller}. 
Observe that recurrence corresponds to local survival if we assume that the
process can not die out globally. Our main result is the criterion
for global survival, Theorem~\ref{main}. This criterion is
given in terms of Lyapunov exponents of an infinite product of
$2\times 2$ random matrices. The main idea of the proof is to
construct an embedded Galton-Watson process in an ergodic random
environment that survives if and only if the BRWRE survives
globally.

This model can be viewed from a different angle. Interpret the
position of a particle as its type. Hence, e.g.\ a particle of type~$0$ 
may produce offspring particles of types $-1,0$, and~$1$. 
%N
Hence, our model is a particular case of a multi-type Galton-Watson process with infinitely many types.
Consider a multi-type Galton-Watson process.
The types are indexed by some set~$I$ that may be finite or countably infinite. 
Then, each type of individual $x\in I$ produces offspring according to
$\mu_x=\mu_x(k_i:~i\in I)$. This means that a particle of
type $x$ has $k_i$ offspring of type $i$ with
probability $\mu_x(k_i)$, $i\in I$.  Assume
irreducibility of the process, i.e., any type of particle may,
after some generations, have any other type of particle as a
descendant. Let  $M:=(m(x,y))_{x,y\in I}$  be the
first moment matrix, i.e., $m(x,y)$ is the mean number of type~$y$
offspring of one type~$x$ particle. If~$I$ is finite, it is
well-known, cf.\ e.g.\ Chapter~II of~\cite{Ha}, 
that the multi-type Galton-Watson process survives with
positive probability if and only if the the largest eigenvalue~$\rho(M)$
 of the matrix~$M$, the Perron-Frobenius eigenvalue, is
greater than~$1$. In the case when~$I$ is infinite, 
very little is known in general.
%N
One reason for this is the following. In the finite case, 
as we just mentioned, the behavior of
the process is strongly connected with the Perron-Frobenius
eigenvalue of the first moment matrix. In the infinite setting
this matrix becomes an operator that does not necessarily have a
\emph{largest eigenvalue}. If the operator~$M$ is ergodic in the sense of~\cite{V},
 most results carry over from the finite case, e.g.\ see
Chapter~III, Section~10 of~\cite{Ha} 
for general results, \cite{BHL} for a
concrete example, and~\cite{V} for mathematical
background. If the operator is not ergodic nothing
is known in general and  the behavior of the process becomes more
subtle. 
%N
There are some interesting results on an example
of epidemics where a classification is obtained, see~\cite{BHL}, 
and further partial results for specific models
in~\cite{BZ} and~\cite{PS}, but no general explicit criterion is known. Therefore, since the first moment matrix of BRWRE is in general not ergodic, cf. Proposition~5.2.8 of~\cite{Muellerthesis}, 
the classification of BRWRE constitutes a step towards the 
understanding of infinite-type Galton-Watson processes in general.
Our model is connected with the model studied in~\cite{BCGH2, GH2, GH3}
(where the random environment was related only to the branching mechanism). In the
model studied there, the process does not start only with one particle but an
infinite number of particles.
While in~\cite{BCGH2, GH2, GH3} the authors analyse different
growth rates using a variational approach, we give a description of
survival and extinction using properties of branching random
walks and embedded Galton-Watson processes in random environment.

% 
% Our model is similar to the one investigated in~\cite{BCGH2, GH2, GH3} 
% (where the random environment was related only to the branching mechanism)
% and we hope to explore the connections in future work. 
% However, our paper is self-contained and does not rely on the results in \cite{BCGH2, GH2, GH3}.

\section{Formal description of the model and main results}
\label{s_model} We now describe the model, keeping the notations
of~\cite{CP, CP07} whenever possible. Let $\Z_+=\{0,1,2,\ldots\}$
and
 $\A=\{-1, 0, 1\}$.
Define
%N
\[
 \V = \Big\{v=(v_x, x\in \A) : v_x \in \Z_+, \forall  x\in \A \Big\},
\]
and for $v\in\V$ put $|v|=\sum_{x\in \A}v_x$.
Furthermore, let~$\M$ be the set of all probability
measures~$\omega$ on $\V$:
\[
 \M = \Big\{\omega = (\omega(v), v\in\V) : \omega(v)\geq 0
\mbox{ for all }v\in\V,
                     \sum_{v\in\V}\omega(v)=1\Big\}.
\]
%N
Then, suppose that $\omega:=(\omega_x\in\M, x\in\Z)$ is an i.i.d.\ sequence with values in $\M$,
and denote by $\IP, \IE$ the probability and expectation with respect to~$\omega$.

The collection $\omega = (\omega_x, x\in\Z)$ is called \emph{the
environment}. Given the environment~$\omega$, the evolution of the
process is described in the following way: start with one particle
at some fixed site of~$\Z$. At each integer time the particles
branch independently in the following way: for a particle
at site $x\in\Z$, a random element $v=(v_y, y\in \A)$ is chosen
with probability $\omega_x(v)$, and then the particle is
substituted by $v_y$ particles in $x+y$, $y\in \A$.

It is important to note that, in contrast to~\cite{CP, CP07}, the
condition $|v|\geq 1$ is dropped from the definition of~$\V$. This
means that here we allow the possibility that particles can
disappear (i.e., leave no offspring), thus it can happen that the process dies out.

Denote
\[
 \mu^-_x=\sum_{v\in\V}\omega_x(v)v_{-1}, \;
\mu^0_x=\sum_{v\in\V}\omega_x(v)v_{0}, \mbox{ and } \;
\mu^+_x=\sum_{v\in\V}\omega_x(v)v_{1}.
\]
In words:  $\mu^-_x$ is the mean number of offspring sent by a particle from~$x$ to $x-1$,
 $\mu^+_x$ is the mean number of offspring sent by a particle from~$x$ to $x+1$,
and~$\mu^0_x$ is the mean number of offspring which stay at~$x$.

We always assume that the following two conditions hold:

\medskip
\noindent
\textbf{Condition~E.} We have $\IP[\min(\mu^-_0, \mu^+_0)>0]=1$.

\medskip
\noindent
\textbf{Condition~B.} 
There exists $v\in\V$ with $|v|\geq 2$ such that $\IP[\omega_0(v)>0]>0$.

\medskip

Condition~E is a natural ellipticity condition which ensures that
the process is irreducible in the sense that for any $x,y\in\Z$ a
particle from~$x$ can have descendants in~$y$. Condition~B says
that there are places where particles are able to branch.

For the proof of Theorem \ref{main} below, we will need the following stronger condition.
Let $V_1: = \{v\in \V: v_1 \geq 1\}$ and $V_{-1}: = \{v\in \V: v_{-1} \geq 1\}$.

\medskip
\noindent
\textbf{Condition~S.} 
$\IE |\ln(\omega_0(V_1))|< \infty$ and $\IE |\ln(\omega_0(V_{-1}))| < \infty$.

\medskip
Since $\mu_0^-\geq \omega_0(V_{-1})$ and $\mu_0^+\geq \omega_0(V_1)$,
Condition~S implies that $\IE (\ln \mu_0^-)^-$ and $\IE (\ln \mu_0^+)^-$ are finite.

Let us denote by~$\eta_n(y)$  the number of particles in~$y$ at time~$n$.
Define the random variable
\[
 \ZZ_n = \sum_{y\in\Z} \eta_n(y),
\]
i.e., $\ZZ_n$ is the total number of particles at time~$n$.

We denote by $\Po^x, \Eo^x$  the 
probability and expectation for the process starting from~$x$ in the fixed
%N
environment~$\omega$, often denoted as ``quenched'' probability and expectation.
% We use the notation
% $\PA^x[\,\cdot\,] = \IE\,\Po^x[\,\cdot\,]$
% for the annealed law of the branching random walk in random environment, and
% $\EA^x$ for the corresponding expectation.
%N
%Also, sometimes we use the symbols $\Po,\Eo$ without the
%corresponding superscripts when it can create no confusion (e.g.,
%when the starting point of the process is indicated elsewhere). 
%We don't need this!

Now we define the survival regimes.

\begin{df}
\label{global}
 Given $\omega$, we say that there is \emph{global survival}
if
\[
\Po^0[ \ZZ_n\to 0]<1.
\]
\end{df}

\begin{df}
\label{local} Given $\omega$,  we say that there is \emph{local
survival} if
\[
\Po^0[ \eta_n(y)\to 0]<1
\]
for all $y$.
\end{df}

\begin{df}
\label{s-local}
 Given $\omega$,  we say that there is \emph{strong local survival}
if
\[
\Po^0[ \ZZ_n\to 0]=\Po^0[ \eta_n(y)\to 0]<1
\]
for all~$y$.
\end{df}

We say that for a given $\omega$ there is local (respectively global) extinction, if there is no  local (respectively global) survival. 

In principle, in a (properly constructed) deterministic environment the
definition of (strong) local survival may depend on the starting point,
cf.\ Example~1 of~\cite{CP}. Let us show, however, that in i.i.d.\
random environment there is no such dependence and that local survival always implies strong local survival.
\begin{theo}
\label{s_loc-ext} Local survival and strong local survival do 
not depend on the starting point in Definitions \ref{local} and \ref{s-local}.
  Also, either there is strong
local survival for $\IP$-a.a.\ $\omega$, or there is local extinction
for $\IP$-a.a.\ $\omega$.
\end{theo}

Similarly to~\cite{CMP, MP1, MP2, Mueller} we can obtain a simple
and explicit criterion for local extinction. As in the above
references, it turns out that local extinction does not depend
on the environmental law itself, but only on its support.

\begin{df}
%N
\label{transient} We say that the process vanishes on the right
(respectively, on the left) if for any $z\in\Z$, the set
$\{z,z+1,z+2,\ldots\}$ (respectively, $\{z,z-1,z-2,\ldots\}$) is
visited only finitely many times a.s.
\end{df}

The criterion for local extinction is then given by
\begin{theo}
 \label{loc_ext}
There is local extinction iff there exists $\lambda>0$ such that
\begin{equation}
 \label{lambda}
\mu_0^-\lambda^{-1} + \mu_0^0+\mu_0^+\lambda\le 1
\end{equation}
for $\IP$-a.a.\ $\omega$. Moreover, if $\lambda>1$, then the process
vanishes on the right, and if $\lambda<1$, then the process 
vanishes on the left.
\end{theo}

\begin{rmk}
If~\eqref{lambda} holds for $\IP$-a.a.\ $\omega$ with $\lambda=1$, then
there is global extinction a.s.\
(so that the process vanishes in both directions).
This is easy to see:
%N
$\mu_0^- + \mu_0^0+\mu_0^+ \le 1$ implies that the
mean offspring in all sites is less than or equal to~$1$, and so
the total number of particles in the process is a (nonnegative) supermartingale.
This supermartingale converges a.s.\ to some limit, and it is straightforward to
obtain (using Conditions~B,~E, and the fact that the environment is i.i.d.)
that this limit can only be~$0$.
\end{rmk}

\begin{rmk}
 By Theorem~\ref{loc_ext}, local extinction implies 
that $\IP[\mu_0^0<1]=1$.
Then, particles cannot accumulate in any site without help from outside.
Using Condition~E we obtain that
if $\IP\big[\omega_0\big((0,0,0)\big)>0\big]>0$, then it is \emph{possible} that the process dies out,
i.e., $\Po^0[\ZZ_n\to 0]>0$.
\end{rmk}

Now, the goal is to obtain a criterion for global extinction in
the case when the process becomes locally extinct. To this end we introduce the following  matrices.
For $k\in \{1, 2, 3, \ldots\}$, denote
\begin{equation}
\label{def_A}
 A_k=\left( \begin{array}{cc}
         \frac{1-\mu_k^0}{\mu_k^+}
             &-\frac{\mu_k^-}{\mu_k^+}\\
          1  &0
        \end{array}
\right),\mbox{ and }\tilde A_k=\left( \begin{array}{cc}
         \frac{1-\mu_k^0}{\mu_k^-}
             &-\frac{\mu_k^+}{\mu_k^-}\\
          1  &0
        \end{array}
\right).
\end{equation}
So, $A_1, A_2, A_3, \ldots$ and $\tilde A_1, \tilde A_2, \tilde A_3, \ldots$ 
are two sequences of i.i.d.\ random matrices.
Denote by~$\gamma_1$ the top Lyapunov exponent associated with the sequence $\{A_k\}$, i.e.,
\[
\gamma_1=\lim_{n\to \infty} \frac{1}{n} \IE(\ln \|A_n\cdots A_1\|),
\]
where $\|\cdot\|$ is any matrix norm (this limit exists provided that 
$\IE\ln^+\|A_1\|$ is finite, cf.\ e.g.\ Section~I.2 of~\cite{BG}). 
Analogously, let $\tilde\gamma_1$ be the top Lyapunov exponent of the sequence $\{\tilde A_k\}$.
%%% QUESTION: REFERENCE WHY THIS LIMIT EXISTS?
The criterion for  global survival is then given by
\begin{theo}
 \label{main}
Suppose that Condition~S holds. Assume also that 
there is local extinction, so, by Theorem~\ref{loc_ext}, there is $\lambda>0$
such that~\eqref{lambda} holds $\IP$-a.s. Then,
\begin{itemize}
\item if there is some $\lambda>1$ such that 
$\mu_0^-\lambda^{-1} + \mu_0^0+\mu_0^+\lambda\le 1$ $\IP$-a.s., then there is global survival iff
\begin{equation}
 \label{condgloballeft}
\gamma_1< \IE\ln\Big(\frac{\mu_0^-}{\mu_0^+}\Big);
\end{equation}
\item if there is some $\lambda<1$ such that $\mu_0^-\lambda^{-1} + \mu_0^0+\mu_0^+\lambda\le 1$ 
$\IP$-a.s., then there is global survival iff
\begin{equation}
 \label{condglobalright}
\tilde\gamma_1< \IE\ln\Big(\frac{\mu_0^+}{\mu_0^-}\Big).
\end{equation}
\end{itemize}
\end{theo}

\section{Proofs}
In this section we prove Theorems~\ref{s_loc-ext}, \ref{loc_ext}, and~\ref{main}.

\subsection{Proof of Theorem \ref{s_loc-ext}.}
%N
An important object is the first moment matrix
$M_\omega=(m_\omega(x,y))_{x,y\in\Z}$ of the process which is
defined by
\[
m_\omega(x,x-1)=\mu_x^-, \quad m_\omega(x,x)=\mu_x^0, \quad m_\omega(x,x+1)=\mu_x^+,
\]
%N
and $m_\omega(x,y)=0$ for $y\notin \{x-1,x,x+1\}$. 
Let $M^n_\omega=(m_\omega^{(n)}(x,y))_{x,y\in \Z}$ denote the
%N
$n$-fold convolution of~$M_\omega$; in other words,
$m^{(n)}(x,y)=\Eo^x[\eta_n(y)]$. Due to Condition~E, the matrix $M_\omega$ is
irreducible. We have, by a supermultiplicativity argument, that
\begin{equation}
\label{spec}
\rho(M_\omega):=\limsup_{n\to\infty}
\left(m_\omega^{(n)}(x,y)\right)^{1/n}
\end{equation} 
does not depend on~$x$ and~$y$ (cf.\ e.g.~\cite{GM}).

Due to the irreducibility Condition~E, we obtain for all $x,z\in\Z$ that 
%N
$\Po^x[z \mbox{ is visited by some particle}]= \Po^x[\eta_n(z) > 0 \mbox{ for some }n ]>0$. Since for $y\in \Z$
\[
\Po^x[\eta_n(y)\not\to 0]\geq \Po^x[z \mbox{ is visited by some
particle}]\times \Po^z[\eta_n(y)\not\to 0],
\] 
%N
local survival does not depend on the choice of the starting point in its definition. To see that the same holds true for
strong local survival, observe first that for all $x,y\in\Z$
\begin{equation}
\label{stronglocal}
\Po^x[\ZZ_n\to 0]\leq \Po^x[\eta_n(y)\to 0].
\end{equation}
We denote by $\eta_n=(\eta_n(x))_{x\in \Z}$ 
the ``global'' configuration of particles at time~$n$. 
Let $\Xi_n$ be the set of all possible particle configurations at time~$n$;
observe that, since we start with one particle at~$0$, this 
set is finite or countably infinite. Now, assume that
\[
\Po^0[ \ZZ_n\to 0]=\Po^0[ \eta_n(y)\to 0].
\]
 Conditioning on the first time step we obtain
\begin{align*}
\Po^0[\ZZ_n\to 0]&= \sum_{\eta\in \Xi_1} \Po^0[\eta_1=\eta] \Po^0[\ZZ_n\to 0 \mid \eta_1=\eta]\\ 
&= \sum_{\eta\in \Xi_1} \Po^0[\eta_1=\eta] \prod_{x\in\{-1,0,1\}} 
 \big(\Po^x[\ZZ_n\to 0]\big)^{\eta(x)}
\end{align*} 
and
\begin{align*}
\Po^0[\eta_n(y)\to 0]&= \sum_{\eta\in \Xi_1} \Po^0[\eta_1=\eta] \Po^0[\eta_n(y)\to 0\mid \eta_1=\eta]\\ 
&= \sum_{\eta\in \Xi_1} \Po^0[\eta_1=\eta] \prod_{x\in\{-1,0,1\}} 
 \big(\Po^x[\eta_n(y)\to 0]\big)^{\eta(x)}
\end{align*}
Therefore, using \eqref{stronglocal},
\[
\Po^x[ \ZZ_n\to 0]=\Po^x[ \eta_n(y)\to 0]
\] 
for  $x\in\{-1,1\}$. Now, the statement follows from an induction argument due to
Condition~E.

The remaining part of the proof splits into three steps (a similar reasoning can be found 
in~\cite{GM}).

\medskip
\noindent{Step 1:} \emph{Local survival is equivalent to
$\rho(M_\omega)>1$.} Let $\rho(M_\omega)>1$ and use the following well-known
approximation property  of the spectral radius,
namely
\begin{equation}\label{approx}
 \rho(M_\omega)=\sup_{|F|<\infty} \rho(M_{\omega,F}).
\end{equation}
Here $M_{\omega,F}$ is the finite matrix over the set $F$ defined
as $m_{\omega,F}(x,y)=m_\omega (x,y)$ for all $x,y\in F$. Due to
\eqref{approx} there exists a finite set $F$ such that
$\rho(M_{\omega, F})>1$.
Since $\rho(M_{\omega,F})\leq\rho(M_{\omega, G})$ 
for $F\subseteq G$, we can choose $F$ to be connected and such that
$0\in F$. Observe that $M_{\omega, F}$ is the first moment matrix
of the multi-type Galton-Watson process that lives on~$F$. This
process can also be interpreted as the embedded process where
particles live only on the set~$F$ and die if they leave this set.
Since $\rho(M_{\omega, F})>1$ this embedded process is
supercritical. This implies the local survival of the BRWRE.

Now, assume local survival of the process. We proceed by constructing an embedded
Galton-Watson process counting the number of particles in the origin. 
Let the particles that are the first particles in their ancestry 
line (of the BRWRE) to return to~$0$ form the first generation of the new process. 
The process is defined inductively: the $i$-th generation consists of particles 
that are the $i$-th particle in their ancestry line to return to~$0$. 
Denote by~$\psi_i$ the size of $i$-th generation.
Observe that $\psi_i\in\N\cup\{\infty\}$ is a Galton-Watson process 
with mean $\Eo \psi_1 > 1$ (in fact, one can even show
that $\Eo \psi_1 =\infty$)
 since the process survives locally. 
Now, we define an embedded process of the above Galton-Watson process,
which is formed by particles that do not go too far away from the origin. 
Let the restricted first generation consist of particles that
 are the first particles in their ancestry 
line to return to~$0$ before time~$N$.
Inductively, the restricted $i$-th generation is formed by the particles
having an ancestor in the restricted $(i-1)$th generation
 and being the first in their 
ancestry line of this ancestor to return to~$0$ in at most~$N$ time steps. 
Let $\psi_i^{(N)}$ be the size of the restricted $i$-th generation and let us
choose~$N$ such that $\Eo \psi_1^{(N)}>1$. 
Setting $\B_N:=[-N,\ldots,N]$ we obtain $\rho(M_\omega)>\rho(M_{\omega,\B_N})>1$.

\medskip
\noindent{Step 2:} \emph{Either there is local survival for
$\IP$-a.a.\ $\omega$, or there is local extinction for
 $\IP$-a.a.\ $\omega$.} The spectral radius $\rho(M_\omega)$ is deterministic. To see this
observe that $\rho(M_\omega)=\limsup_n (m^{(n)}(x,y))^{1/n}$ does
%N
not depend on~$x$ and~$y$ and is
constant $\IP$-a.s.\ by ergodicity of the environment as discussed in~\cite{Mueller}.

\medskip
\noindent{Step 3:} \emph{Local survival implies strong local
survival.} We assume local survival. Recall that a set $F$ with
$\rho(M_{\omega, F})>1$ gives rise to a supercritical multi-type
Galton-Watson process. In analogy to \cite{CP} we call these sets
\emph{(recurrent) seeds}. We make the following observation that
is obvious if $\IP$ is discrete and easy to check otherwise: There
%N
exists some $N\in\N$ and some $\eps>0$
\begin{equation}
\label{eq:goodseed}
 \IP\left[ \PoA{\B_N}^x[\ZZ_n\not\to 0]>\eps\right]>0,
\end{equation}
here $\PoA{\B_N}^x$ denotes the probability measure of the
%N
embedded process that starts in $x\in \B_N$ and lives on~$\B_N$ (i.e. particles leaving $\B_N$ die).
We proceed similarly to the proof of 
Lemma~2.6 of~\cite{CP} and partition~$\Z$ into translates of~$\B_N$. 
% Since the environment is i.i.d., we can construct the
% BRWRE by choosing randomly the environment in a translate of~$\B_N$ 
% at the first moment when the process enters this set.
 Then,
by the Borel-Cantelli Lemma, infinitely many of the translates of~$\B_N$
contain a seed with survival probability at least~$\eps$, cf.~\eqref{eq:goodseed}. 
Now, if the process survives globally infinitely many such seeds will be
visited and it is straightforward to construct an independent
sequence of embedded supercritical multi-type Galton-Watson
processes whose survival probability is greater than~$\eps$.
Eventually, one of those processes survives and strong local
survival follows. \qed

\subsection{Proof of Theorem \ref{loc_ext}.}
First let us observe that  Lemma 3.5 in \cite{Mueller} implies that
%N It is R not R^d here!
\begin{equation}
\rho(M_\omega)=\sup_{\mu\in \hat\K} \inf_{\lambda\in\R} \Big(\lambda
\mu^+ + \mu^0 +\lambda^{-1} \mu^-\Big),
\end{equation} 
where $\mu=(\mu^+,\mu^0,\mu^-)$ and $\hat\K$ is the convex hull of the support of the
one-dimensional marginal of~$\IP$. Observing that~$\sup$ and~$\inf$ 
are attained with say~$\hat \mu$ and~$\hat\lambda$, we
obtain that $\hat\lambda\hat
\mu^+ + \hat \mu ^0 +\hat\lambda^{-1}\hat \mu^- \le 1$ implies~\eqref{lambda}. 
Clearly, \eqref{lambda} implies that
$\inf_{\lambda\in\R} \sup_{\mu\in \hat\K} (\lambda \mu^+ + \mu^0
+\lambda^{-1} \mu^-) \le 1 $. Since, by a
minimax argument, 
one can exchange~$\inf$ and~$\sup$, we obtain that $\rho(M_\omega)\le 1$ a.s.

Now, let us suppose that $\lambda>1$ and let us prove that
the process vanishes on the right. Note that, by~\eqref{lambda}, the function
\[
h(n)=\sum_{x\in\Z}\eta_n(x)\lambda^x
\]
is a positive supermartingale for $\IP$-a.a.\ $\omega$ (see, for example, the
proof of Theorem~1.6  in~\cite{CP}). Therefore, as $n\to\infty$, it converges a.s.\ 
to some random variable~$h_\infty$. Using Fatou's Lemma, 
for the process starting at the origin we
obtain that, $\IP$-a.s.,
\[
 \Eo h_\infty\le \Eo h(0)= 1
\]
On the event that the set $\{1,2,3,\ldots\}$ is visited infinitely often 
we have that every $k\ge 1$ is visited at least once. 
Hence, $\limsup h(n)\ge \lambda^k$ for all~$k$. 
 Since this contradicts the fact that~$h(n)$ converges to a finite random variable, 
we obtain that $\{1,2,3,\ldots\}$ is only visited finitely often a.s.
Using the irreducibility we obtain that the set $\{z,z+1,z+2,\ldots\}$
is visited finitely often a.s.\ for any $z\in\Z$. The case $\lambda<1$
can be treated analogously.
 \qed

\subsection{Proof of Theorem~\ref{main}.}
Let us assume  that there exists $\lambda>1$ which satisfies~\eqref{lambda} for $\IP$-a.a.\
%N
$\omega$ (so, the process vanishes on the right). The proof for the case $\lambda < 1$
follows then by exchanging~$\lambda$ and~$\lambda^{-1}$ and~$\mu^-$ and~$\mu^+$ 
(that is, using $\tilde A_k$ instead of $A_k$).

Since the matrices $\{A_k\}$ are not nonnegative,
 we introduce the following sequences of nonnegative matrices. 
For $k\in \{1, 2, 3, \ldots\}$, denote
\begin{equation}
\label{def_Ak}
 A_k^{(\lambda)}=\left( \begin{array}{cc}
         \frac{\mu_k^-}{\lambda^2\mu_k^+}
             &\frac{1-\mu_k^0-\lambda^{-1}\mu_k^--\lambda\mu_k^+}{\lambda \mu_k^+}\\
          \frac{\mu_k^-}{\lambda^2\mu_k^+}
                &1+\frac{1-\mu_k^0-\lambda^{-1}\mu_k^--\lambda\mu_k^+}{\lambda \mu_k^+}
        \end{array}
\right).
\end{equation}
That is, $A_1^{(\lambda)}, A_2^{(\lambda)}, A_3^{(\lambda)}, \ldots$ is a sequence 
of i.i.d.\ random matrices, which are nonnegative by~\eqref{lambda}.
Denote by~$\gamma_1^{(\lambda)}$ the top Lyapunov exponent associated with the sequence $\{A_k^{(\lambda)}\}$
and by~$\gamma_2^{(\lambda)}$ the second Lyapunov exponent. It holds that (cf., for example,
Corollary~1.3 of~\cite{L})
\begin{equation}
 \label{sum_exp}
\gamma_1^{(\lambda)}+\gamma_2^{(\lambda)}=\IE\Big[\ln \det A_1^{(\lambda)}\Big]=\IE\ln\Big(\frac{\mu_k^-}{\lambda^2\mu_k^+}\Big)
=\IE\ln\Big(\frac{\mu_k^-}{\mu_k^+}\Big) -2\ln \lambda.
\end{equation}

The proof splits into two parts. First, we prove that there is global survival iff
\begin{equation} \label{cond_gamma}
\gamma_1^{(\lambda)}< \IE\ln\Big(\frac{\mu_0^-}{\mu_0^+}\Big) -\ln \lambda.
\end{equation}
We conclude then by comparing the Lyapunov spectra of $\{A_k\}$ and $\{A_k^{(\lambda)}\}$.

We will consider two modifications of our BRWRE. The first modification is the following.
Start the original BRWRE with one particle at $0$. When a particle
hits $-1$, it is frozen and remains at $-1$ until all the existing
particles hit $-1$ (this will happen in finite time, as our
process vanishes on the right). Let $\xi_1$ be the total number
of  frozen particles at $-1$. Then, release the frozen particles,
let them perform a BRW in random environment $\omega$, and freeze all
particles that hit $-2$. When all the existing particles are
frozen at $-2$, let $\xi_2$ be the number of particles at $-2$. We
repeat the above construction in this way to obtain a branching process $\{\xi_n\}_{n=1,2,\ldots}$ in stationary ergodic random environment. By Theorem 5.5 and Corollary 6.3 of~\cite{T} (taking into account Condition~S),
the above process survives with positive probability iff
$\IE\ln\Eo\xi_1>0$. 
%The two processes $\{\xi_n\}$ and $\{\ZZ_n\}$ have the same trajectories up to a time change 
%{\tt Nina: this is not really true. I suggest:}
But survival of the process
$\{\xi_n\}$ is equivalent to survival of our original
process $\{\ZZ_n\}$.  

We are going to calculate $\IE\ln\Eo\xi_1$ by means of constructing another
sequence of random variables $\{\zeta_n\}$ in such a way that $\IE\ln\Eo\xi_1=\IE\ln\Eo\zeta_1$.

Now, our second modification of the original BRWRE is defined in
the following way. 
%N
We start with one particle in $k \geq 0$. When a particle hits~$0$, it is frozen and
remains at~$0$ forever, i.e., we modify the environment by putting $\omega'_0(v')=1$, where
$v'=(0,1,0)$. Denote by~$\zeta_k$ the total number of frozen
particles at~$0$ starting with one particle at~$k$.  As the
environment~$\omega$ is stationary, $\zeta_k$ and~$\xi_k$ have the
same annealed law, $k=1,2,3,\ldots$.
% $\Eo\zeta_k$ and
% $\Eo\xi_k$ have the same distribution, $k=1,2\ldots$
%  As the environment $\omega$ is stationary, we have
% $\IE\ln\Eo\xi_1=\IE\ln\Eo\zeta_1$.

Denote $f(k)=\Eo^k \zeta_k$, $k=0, 1, 2, \ldots$. Note that $f(0)=1$,
%N the equation is not recurrent! :-)
and for $k\geq 1$ we can write the recursive equation
\begin{equation}
\label{rec_f}
f(k)=\mu_k^- f(k-1)+\mu_k^0 f(k)+\mu_k^+ f(k+1).
\end{equation}
Let $g(k)=\lambda^{-k} f(k)$. Then,~\eqref{rec_f} implies that
\begin{equation}
\label{rec_g}
\lambda^k g(k)=\lambda^{k-1}\mu_k^- g(k-1)+\lambda^k\mu_k^0 g(k)+\lambda^{k+1}\mu_k^+ g(k+1).
\end{equation}
Denote $\Delta(k)=g(k)-g(k-1)$. Observe that~\eqref{rec_g} can be rewritten
as
\begin{align*}
 \left( \begin{array}{c}
         \Delta(k+1)\\
         g(k+1)
        \end{array}
\right)
  =A_k^{(\lambda)}  \left( \begin{array}{c}
         \Delta(k)\\
         g(k)
        \end{array}
\right),
\end{align*}
where $A_k^{(\lambda)}$ is the matrix defined in~\eqref{def_Ak}.
Recall that, by~\eqref{lambda}, the matrix~$A_k^{(\lambda)}$ is nonnegative.

In fact, to define Lyapunov exponents and use the classical results about them, we need
$\ln \|A_1^{(\lambda)}\|$ and $\ln \|(A_1^{(\lambda)})^{-1}\|$ to be integrable.
It  is straightforward to check
that this is the case iff $\IE (\ln \mu_k^-)^-$ and  $\IE (\ln \mu_k^+)^-$ are finite,
which is a consequence of Condition~S.
% As all norms in matrix space are equivalent, take
% $\|A\|=\sum a_{ij}$. Then,
% \begin{align*}
%  \|A_k\|&=1+\frac{2\mu_k^-}{\lambda^2 \mu_k^+}
%      +\frac{ 2(1-\lambda^{-1}\mu_k^--\mu_k^0-\lambda\mu_k^+)}{\lambda\mu_k^+ },\\
%  \|A_k^{-1}\|&=2+\frac{ 1-\lambda^{-1}\mu_k^--\mu_k^0}{\lambda^{-1}\mu_k^-}
%      +\frac{ 1-\lambda^{-1}\mu_k^--\mu_k^0-\lambda\mu_k^+}{\lambda^{-1}\mu_k^- }.
% \end{align*}
% As $0<1-\lambda^{-1}\mu_k^--\mu_k^0-\lambda\mu_k^+<1$, the bad situation is only when
% $\mu_k^-$ or $\mu_k^+$ are small. So, we need $\ln(\mu_k^-)^{-1}$,  $\ln(\mu_k^+)^{-1}$
% to be integrable at $0$.

Denote by $H_\omega\subset \R^2$ the random one-dimensional subspace
of $\R^2$ associated with $\gamma_2^{(\lambda)}$.
 So, for all $e\in \R^2\setminus H_\omega$, we have
\begin{equation}
\label{g1}
\lim_{n\to \infty} \frac{1}{n} \ln \|A_n^{(\lambda)}\cdots A_1^{(\lambda)} e\|=\gamma_1^{(\lambda)},
\end{equation}
and for all $e'\in H_\omega\setminus\{0\}$ we have
\begin{equation}
\label{g2}
\lim_{n\to \infty} \frac{1}{n} \ln \|A_n^{(\lambda)}\cdots A_1^{(\lambda)} e'\|=\gamma_2^{(\lambda)}
\end{equation}
(cf., for example, Theorem 3.1 of~\cite{L}).
If $\mu_0^-\lambda^{-1} + \mu_0^0+\mu_0^+\lambda< 1$ with positive probability, then
Corollary 2 of~\cite{H} implies that $\gamma_2^{(\lambda)}<\gamma_1^{(\lambda)}$.
If, on the other hand, $\mu_0^-\lambda^{-1} + \mu_0^0+\mu_0^+\lambda= 1$ $\IP$-a.s.,
 after some elementary computations it is straightforward to obtain that
 $\{\gamma_1^{(\lambda)}, \gamma_2^{(\lambda)}\}=\{0, \IE\ln(\mu_k^-/(\lambda^2\mu_k^+))\}$ and thus we also have
 $\gamma_2^{(\lambda)}<\gamma_1^{(\lambda)}$ (or $\gamma_1^{(\lambda)}=\gamma_2^{(\lambda)}=0$;
we treat this case later).

So, suppose that $\gamma_2^{(\lambda)}<\gamma_1^{(\lambda)}$. 
%N
Now, our goal is to prove that
$(\Delta(1), g(1))\in H_\omega$.  We argue by contradiction.
Observe that there exists some vector $e_\omega\ge (1,1)$ such that 
 $e_\omega\notin H_\omega$ and $e_\omega \cdot (\Delta(1), g(1)) \neq 0$.
%N We need this, otherwise there is no $c_\omega$.
Suppose that $\phi_\omega:=(\Delta(1), g(1))\notin H_\omega$.
Then, there exists a number $c_\omega\ne 0$ such that $u_\omega:=\phi_\omega +c_\omega e_\omega\in
H_\omega$. Now, let us write
\begin{align}
 \label{g1_wins}
 \left( \begin{array}{c}
         \Delta(k+1)\\
         g(k+1)
        \end{array}
\right)
  =\Big(\prod_{i=1}^k A_i^{(\lambda)}\Big) \phi_\omega=
    \Big(\prod_{i=1}^k A_i^{(\lambda)}\Big)u_\omega
  -c_\omega\Big(\prod_{i=1}^k A_i^{(\lambda)}\Big)e_\omega.
\end{align}
As $\gamma_1^{(\lambda)}>\gamma_2^{(\lambda)}$, using~\eqref{g1} and~\eqref{g2}, 
and the uniform positiveness of $e_\omega$, we see that for all $k$ large enough
\begin{equation}
\label{sign}
\sgn \Delta(k)=\sgn g(k)=-\sgn (c_\omega).
\end{equation}

On the other hand, let us show that for all $k$  we have $\Delta(k)\le 0$.
Then, since $g(k)>0$ by definition, we obtain a contradiction with~\eqref{sign}.

Indeed, we have
\[
\Delta(k)=g(k)-g(k-1)=\lambda^{-k} \Eo^k \zeta_k -\lambda^{-(k-1)} \Eo^{k-1} \zeta_{k-1}.
\]
Thus, we need to show that
\begin{equation}
\label{zetas}
\Eo^k \zeta_k \le\lambda \Eo^{k-1} \zeta_{k-1}.
\end{equation}
Note that
\[
 \Eo^k \zeta_k =\Eo^{k-1} \zeta_{k-1} \; \Eo \hat\zeta_k,
\]
where $\hat\zeta_k$ is a random variable defined as follows:
start the process with one particle at $k$ and  freeze all particles that reach $k-1$; then,
$\hat\zeta_k$ is the number of frozen particles at $k-1$.
Observe that, similarly to the proof of Theorem~\ref{loc_ext}, the function
\[
h(\eta_n)=\sum_{x\in\Z}\eta_n(x)\lambda^x
\]
is still a supermartingale for $\IP$-a.a.\ $\omega$ for this process as well.  
So, suppose that we start from
one particle in~$1$ and  freeze all particles that reach~$0$ and
let~$\tau$ be the moment when all particles are frozen. As we
assumed that the cloud of particles vanishes on the right,
$\tau$ is finite a.s. Then, using  Fatou's Lemma, we
obtain that, $\IP$-a.s.,
\[
 \Eo \zeta_1=\Eo \hat\zeta_1=\Eo h(\eta_\tau)\le h(\eta_0)=\lambda.
\]
By stationarity, $\Eo\hat\zeta_k\le \lambda$ for all
$k=1,2,3,\ldots,$ and this shows~\eqref{zetas}.

Hence,
$(\Delta(1), g(1))\in H_\omega$ and
\begin{align}
\label{gA}
 \left( \begin{array}{c}
         \Delta(k+1)\\
         g(k+1)
        \end{array}
\right)
  =\Big(\prod_{i=1}^k A_i^{(\lambda)}\Big) \left( \begin{array}{c}
         \Delta(1)\\
         g(1)
        \end{array}\right).
\end{align}
Let  $\|\cdot \|_1$ be the $L_1$-norm in $\R^2$. Then,
\[
\|(\Delta(k+1), g(k+1))\|_1=|g(k+1)-g(k)|+g(k+1)=g(k),
\]
as $\Delta(k+1)\le 0$. Thus,~\eqref{gA} and~\eqref{g2} imply that
\begin{equation}
\label{asimp_g}
\lim_{k\to \infty} \frac{1}{k} \ln g(k)=\gamma_2^{(\lambda)}
\end{equation}
and so
\begin{equation}
\label{asimp_f}
\lim_{k\to \infty} \frac{1}{k} \ln f(k)
%=\ln \lambda\lim_{k\to \infty} \frac{1}{k} \ln g(k)
=\gamma_2^{(\lambda)}+\ln\lambda.
\end{equation}
As mentioned above, it may happen also that $\gamma_1^{(\lambda)}=\gamma_2^{(\lambda)}=0$,
but in this case it is straightforward to obtain that~\eqref{asimp_g} and~\eqref{asimp_f}
hold as well (since the limits in~\eqref{g1} and~\eqref{g2} are both equal to~$0$).

Now, note that, $\Eo^k\zeta_k = \Eo^k \hat\zeta_k\Eo^{k-1} \hat\zeta_{k-1}\ldots\Eo^1 \hat\zeta_1$,
and $\hat\zeta_1,\ldots,\hat\zeta_k$ have the same annealed law as~$\zeta_1$.
Therefore, by the Ergodic Theorem, we have
\[
\lim_{k\to \infty} \frac{1}{k} \ln f(k)
=\IE\ln \Eo\zeta_1, \quad \text{$\IP$-a.s.}
\]
Thus, $\IE\ln \Eo\zeta_1=\ln\lambda+\gamma_2^{(\lambda)}$, and the process
${\ZZ_n}$ survives globally iff
\[
0<\ln\lambda+\gamma_2^{(\lambda)}=\IE\ln\Big(\frac{\mu_k^-}{\mu_k^+}\Big)
  - \ln \lambda - \gamma_1^{(\lambda)},
\]
by~\eqref{sum_exp}.
Now condition \eqref{cond_gamma} follows. It remains to  prove that
$\gamma_1=\gamma_1^{(\lambda)}+\ln \lambda$. 
Observe that \eqref{rec_f} can be written in terms of~$\{A_k\}$:
\begin{align*}
%\label{fA}
 \left( \begin{array}{c}
         f(k+1)\\
         f(k)
        \end{array}
\right)
  = A_k  \left( \begin{array}{c}
         f(k)\\
         f(k-1)
        \end{array}\right).
\end{align*} 
Using this relation, it is straightforward to check that
$ A_k=\lambda B^{-1} A_k^{(\lambda)}  B,$ where
\[
B:=\left( \begin{array}{cc}
         1 & -\lambda \\
        1 & 0
        \end{array}\right).
\] Now the desired statement follows since $\prod_{k=1}^n A_k =\lambda^n B^{-1} \prod_{k=1}^n A_k^{(\lambda)} B$.
\qed

\section*{Acknowledgements}
S.P.\ and M.V.\ are grateful to Fapesp (thematic grant 04/07276--2), 
CNPq (grants 300328/2005--2, 304561/2006--1, 471925/2006--3) for financial support.
S.M. thanks DFG (project MU 2868/1--1) for financial support.
All authors thank CAPES/DAAD (Probral) for support. 
We thank Christian Bartsch and Michael Kochler for pointing out a mistake in a previous version.

\end{document}